\documentclass[12pt,a4paper]{amsart}
\usepackage{amssymb,cite}
\usepackage{graphicx}
\usepackage{dcpic,pictex}
\usepackage[dvipsnames]{xcolor}

\usepackage{listings}
\theoremstyle{plain}

\newtheorem{labeling}{Labeling}
\newtheorem{algorithm}{Algorithm}
\theoremstyle{definition}

\theoremstyle{remark}

\numberwithin{equation}{section}

\numberwithin{table}{section}

\numberwithin{figure}{section}

\begin{document}

\title[The $\mu$-permanent, a new graph labeling, and a known integer sequence]{The $\mu$-permanent,
a new graph labeling, and a known integer sequence}

\author{Milica An\dj eli\'c}
\address{Department of Mathematics, Kuwait University, Safat 13060, Kuwait}
\email{milica@sci.kuniv.edu.kw}

\author{Carlos M. da Fonseca}
\address{Department of Mathematics, Kuwait University, Safat 13060, Kuwait}
\email{carlos@sci.kuniv.edu.kw}

\author{Ant\'onio Pereira}
\address{Departamento de Matem\'atica, Universidade de Aveiro, 3810-193 Aveiro, Portugal }
\email{antoniop@ua.pt}


\date{March 28, 2016}


\begin{abstract}
Let $A=(a_{ij})$ be an $n$-by-$n$ matrix. For any real number $\mu$,
we define the polynomial
$$P_\mu(A)=\sum_{\sigma\in S_n}
a_{1\sigma(1)}\cdots a_{n\sigma(n)}\,\mu^{\ell(\sigma)}\; ,$$ as the
$\mu$-permanent of $A$, where $\ell(\sigma)$ is the number of
inversions of the permutation $\sigma$ in the symmetric group $S_n$.
In this note, motivated by this notion, we discuss a new graph
labeling for trees whose matrices satisfy certain $\mu$-permanental
identities.  We relate the number of labelings of a path with a
known integer sequence. Several examples are provided.
\end{abstract}

\maketitle

\section{Introduction}

Given an $n\times n$ matrix $A=(a_{ij})$ and a real number $\mu$, we
define the $\mu$-permanent of $A$ as the polynomial
\begin{equation}\label{mup}
P_\mu(A)=\sum_{\sigma\in S_n} \left(\prod_{i=1}^{n}
a_{i\sigma(i)}\right)\mu^{\ell(\sigma)}\; ,
\end{equation}
where $\ell(\sigma)$ is the number of inversions of the permutation
$\sigma$ in the symmetric group $S_n$ of degree $n$, i.e., the
number of interchanges of consecutive elements necessary to arrange
$\sigma$ in its natural order \cite[p.1]{Muir} or, equivalently,
$$\ell(\sigma)=\#\{(i,j)\in \{1,\ldots,n\}^2\, |\, i<j \mbox{ and }
\sigma(i)>\sigma(j)\}\, .$$ For example, we have
$$
  P_\mu
\left(
  \begin{array}{cc}
    a_{11} & a_{12}   \\
    a_{12} & a_{22}
  \end{array}
\right)
   = a_{11}a_{22}+a_{12}^2\mu
$$
and
$$
  P_\mu
\left(
  \begin{array}{ccc}
    a_{11} & a_{12} & a_{13} \\
    a_{12} & a_{22} & a_{23} \\
    a_{13} & a_{23} & a_{33} \\
  \end{array}
\right)
   = a_{11}a_{22}a_{33}+a_{11}a_{23}^2\,\mu+a_{12}^2a_{33}\,\mu+2 a_{12}a_{23}a_{13}\, \mu^2+a_{13}^2a_{22}\,\mu^3\,
   .
$$

The $\mu$-permanent of a square matrix is a natural extension of the
determinant (setting $\mu=-1$) and the permanent (setting $\mu=1$)
which is in fact quite hard to compute \cite[p.190]{BR}. In
addition, making $\mu=0$, we get the product of the main diagonal
entries of the matrix.

This concept was introduced independently and almost simultaneously
by different authors by a quarter of a century ago under different
names in matrix theory and, surprisingly, in Grassmann algebras and
quantum groups: $q$-permanent is just one of them
\cite{B,BL,BR,HH1992,L,L2,NYM1993,NYM1989,T1993,T1991,Y1991}. Here
we adopt one of the possible ways to call this matricial function
(cf. also \cite{dF201x,dF2,dF2005}).

Listing \ref{list1} presents a Mathematica \cite{MM} routine for
computing the $\mu$-permanent of a square matrix.

\begin{lstlisting}[float=htb, label=list1,caption=Routine to compute the $\mu$-permanent]

inversionList[s_] := Module[{i, inverse = Ordering[s]},
  Table[Length[Select[Take[s, inverse[[i]]], (# > i)&]],
    {i,Length[s] - 1}]]

inversions[s_] := Apply[Plus, inversionList[s]]

permanent[A_, mu_] := Module[{n = Length[A]},
  Sum[Product[A[[i, s[[i]]]], {i, n}]* mu^inversions[s],
    {s, Permutations[Range[n]]}]]
\end{lstlisting}

We believe that this code will be particularly useful for further
developments on the properties of the $\mu$-permanent.

It is clear that in general, under similarity, the $\mu$-permanent
does not keep the same value, i.e., the polynomial $P_\mu(A)$ is not
necessarily the same as $P_\mu(BAB^{-1})$, for $B$ nonsingular. In
particular, for permutation similarity. This means that
interchanging rows and columns of the same indexes leads to possible
different $\mu$-permanents. Since interchanging rows and columns
does not change the underlying graph of the matrices involved, but
the labeling of the vertices, we conclude that the $\mu$-permanent
of a graph depends on its labeling.

After a first attempt to extend monotonic properties of the
$\mu$-permanent of Jacobi positive definite matrices to more general
acyclic matrices \cite{dF2005}, it has recently been noticed a
particular labeling for which the previous properties were indeed
satisfied \cite{dF201x}.

In this note we aim to discuss this new labeling for trees, counting
them for paths. Incidentally, this process will lead to a new
interpretation for a well-known integer sequence.

\section{A new graph labeling}

Given a symmetric matrix $A$, the graph of $A$ is defined by the
zero-nonzero off-main diagonal pattern of $A$. In general, the
vertex labeling is not much discussed in matrix theory since most of
the results involve the spectra of matrices, which do not change by
such labelings. For example, the underlying graph of a tridiagonal
matrix is a path with the vertices ordered successively $1,2,\ldots,
n$ and edges joining consecutive vertices $i$ and $i+1$:

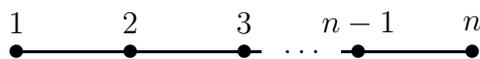
\begin{center}
\setlength{\unitlength}{0.75mm} \thicklines
\begin{picture}(80,15)

\put(0,0){\circle*{2.5}} \put(20,0){\circle*{2.5}}
\put(40,0){\circle*{2.5}} \put(60,0){\circle*{2.5}}
\put(80,0){\circle*{2.5}}

\put(0,0){\line(1,0){43}}

\put(50,0){\makebox(0,0){$\ldots$}}

\put(57,0){\line(1,0){23}}

\put(0,5){\makebox(0,0){$1$}} \put(20,5){\makebox(0,0){$2$}}
\put(40,5){\makebox(0,0){$3$}} \put(60,5){\makebox(0,0){$n-1$}}
\put(80,5){\makebox(0,0){$n$}}

\end{picture}
\end{center}

\hspace{1em}

However, as we mentioned in the introduction, this is not the case for the $\mu$-permanent of a square matrix. For instance, we have
$$
  P_\mu
\left(
  \begin{array}{ccc}
    a_{11} & a_{12} & 0 \\
    a_{12} & a_{22} & a_{23} \\
   0 & a_{23} & a_{33} \\
  \end{array}
\right)
   = a_{11}a_{22}a_{33}+a_{11}a_{23}^2\,\mu+a_{12}^2a_{33}\,\mu\,
   ,
$$
which is a polynomial of degree $1$, and
$$
  P_\mu
\left(
  \begin{array}{ccc}
    a_{11} & a_{12} & a_{13} \\
    a_{12} & a_{22} & 0 \\
    a_{13} & 0 & a_{33} \\
  \end{array}
\right)
   = a_{11}a_{22}a_{33}+a_{12}^2a_{33}\,\mu+a_{13}^2a_{22}\,\mu^3\,
   ,
$$
which in turn has degree $3$. The ``graph'' is the same but the
labeling is not, i.e.,
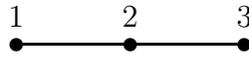
\begin{center}
\setlength{\unitlength}{0.75mm} \thicklines
\begin{picture}(40,10)

\put(0,0){\circle*{2.5}} \put(20,0){\circle*{2.5}}
\put(40,0){\circle*{2.5}}

\put(0,0){\line(1,0){40}}

\put(0,5){\makebox(0,0){$1$}} \put(20,5){\makebox(0,0){$2$}}
\put(40,5){\makebox(0,0){$3$}}
\end{picture}
\end{center}

and

\begin{center}
\setlength{\unitlength}{0.75mm} \thicklines
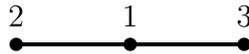
\begin{picture}(40,10)

\put(0,0){\circle*{2.5}} \put(20,0){\circle*{2.5}}
\put(40,0){\circle*{2.5}}

\put(0,0){\line(1,0){40}}

\put(0,5){\makebox(0,0){$2$}} \put(20,5){\makebox(0,0){$1$}}
\put(40,5){\makebox(0,0){$3$}}

\end{picture}
\end{center}

\noindent respectively.

In order to establish several general results for the
$\mu$-permanent, recently in \cite{dF201x}, the second author
introduced the following labeling:

\begin{labeling} Given two disjoint edges $ij$ and $k\ell$, say $i<j$, $k<\ell$,
and $i<k$, then one of the following conditions must be fulfilled:
\begin{itemize}
  \item[(i)] $i<j<k<\ell$
  \item[(ii)] $i<k<\ell<j$.
\end{itemize}
\end{labeling}

To the best of our knowledge, this labeling is new and will be
referred to in the remaining of this paper as a $\mu$-labeling. For
example, we have the following:

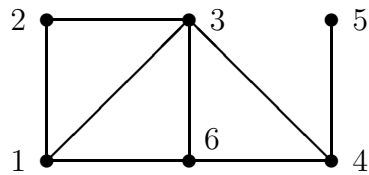
\begin{center}
\setlength{\unitlength}{0.75mm} \thicklines
\begin{picture}(65,35)

\put(5,0){\circle*{2.5}} \put(30,0){\circle*{2.5}}
\put(5,25){\circle*{2.5}} \put(30,25){\circle*{2.5}}
\put(55,25){\circle*{2.5}} \put(55,0){\circle*{2.5}}

\put(5,0){\line(1,1){25}} \put(5,0){\line(1,0){50}}
\put(5,0){\line(0,1){25}} \put(55,0){\line(-1,1){25}}
\put(55,0){\line(0,1){25}}

\put(5,25){\line(1,0){25}} \put(30,0){\line(0,1){25}}

\put(0,0){\makebox(0,0){$1$}} \put(60,0){\makebox(0,0){$4$}}
\put(0,25){\makebox(0,0){$2$}} \put(60,25){\makebox(0,0){$5$}}
\put(34,4){\makebox(0,0){$6$}} \put(35,25){\makebox(0,0){$3$}}

\end{picture}
\end{center}

\hspace{1em}

Interestingly, not all graphs allow such labeling: a complete graph
of more than $3$ vertices is just an example. However, any tree
allows labelings satisfying the conditions described. In what
follows, we present an algorithm to construct one of such labeling:

\begin{algorithm}
Let us consider a tree with a given number of vertices.
\begin{enumerate}
  \item[Step 1.] Choose any vertex from the tree and label it with $1$, which will be the root.
  \item[Step 2.] Take the largest path attached to vertex $1$.
  \item[Step 3.] Label the vertices of this path by
  $2,3,\ldots,k$, where $(i,i+1)$ is an edge, for $i=1,2,\ldots,k-1$.
  \item[Step 4.] Choose the vertex with largest label in the previous with degree greater
  than two, say $\ell$.
  \item[Step 5.] Repeat step 2., replacing $1$ by $\ell$.
  \item[Step 6.] Repeat step 3., labeling the vertices of the path by $k+1,\ldots,k'$.
  \item[Step 7.] Once all vertices of degree more than two were
  considered, restart from 2., choosing the second largest path
  attached to the root and proceed until all vertices were considered.
\end{enumerate}
\end{algorithm}

We remark that, if there is more than one path attached to the root
of the same largest size we choose arbitrarily one of them. Clearly,
this is not the only way to construct such labeling.

As a simple example of the algorithm, we have

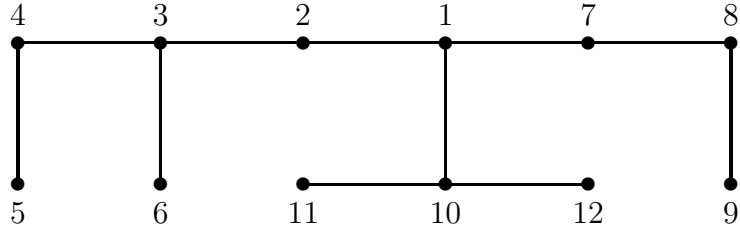
\begin{center}
\setlength{\unitlength}{0.75mm} \thicklines
\begin{picture}(125,40)

\put(0,5){\circle*{2.5}} \put(50,5){\circle*{2.5}}
\put(75,5){\circle*{2.5}} \put(100,5){\circle*{2.5}}
\put(125,5){\circle*{2.5}} \put(25,5){\circle*{2.5}}

\put(0,30){\circle*{2.5}} \put(25,30){\circle*{2.5}}
\put(50,30){\circle*{2.5}} \put(75,30){\circle*{2.5}}
\put(100,30){\circle*{2.5}} \put(125,30){\circle*{2.5}}

\put(0,30){\line(1,0){125}} \put(0,30){\line(0,-1){25}}
\put(25,30){\line(0,-1){25}}
\put(125,30){\line(0,-1){25}}\put(75,30){\line(0,-1){25}}
\put(50,5){\line(1,0){50}}

\put(75,35){\makebox(0,0){$1$}} \put(50,35){\makebox(0,0){$2$}}
\put(25,35){\makebox(0,0){$3$}} \put(0,35){\makebox(0,0){$4$}}
\put(0,0){\makebox(0,0){$5$}} \put(25,0){\makebox(0,0){$6$}}

\put(100,35){\makebox(0,0){$7$}} \put(125,35){\makebox(0,0){$8$}}
\put(125,0){\makebox(0,0){$9$}}

\put(75,0){\makebox(0,0){$10$}} \put(50,0){\makebox(0,0){$11$}}
\put(100,0){\makebox(0,0){$12$}}

\end{picture}
\end{center}

\hspace{1em}

Returning to the $\mu$-permanent, as a consequence, for any $n\times
n$ matrix $A$ whose graph is a tree with the vertices labeled as
described before, one always has
$$P_\mu(A)  =  a_{ii}P_\mu(A_{i})+ \sum_{i\sim j}|a_{ij}|^2
P_\mu(A_{ij}) \, \mu^{\ell(ij)},$$ for any vertex $i$, or
$$
    \frac{d}{d\mu}\, P_\mu(A)  =  \sum_{i\sim j} \ell(ij)
    |a_{ij}|^2 P_\mu(A_{ij})  \mu^{\ell(ij)-1} \; ,$$
with $i<j$, (cf. \cite{dF2005,dF2,dF201x}). Here, $A_S$ is the
matrix obtained from $A$ replacing the rows and columns indexed by
$S$, by zero, except the entries in the main diagonal, which are
1's.

\section{Counting labelings for paths}

In this section we confine our study to paths. Our algorithm
provides the following example for a path of $5$ vertices:

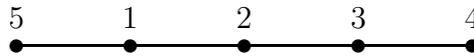
\begin{center}
\setlength{\unitlength}{0.75mm} \thicklines
\begin{picture}(80,15)

\put(0,0){\circle*{2.5}} \put(20,0){\circle*{2.5}}
\put(40,0){\circle*{2.5}} \put(60,0){\circle*{2.5}}
\put(80,0){\circle*{2.5}}

\put(0,0){\line(1,0){80}}

\put(0,5){\makebox(0,0){$5$}} \put(20,5){\makebox(0,0){$1$}}
\put(40,5){\makebox(0,0){$2$}} \put(60,5){\makebox(0,0){$3$}}
\put(80,5){\makebox(0,0){$4$}}

\end{picture}
\end{center}

\hspace{1em}

The following labeling is also a possibility

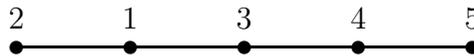
\begin{center}
\setlength{\unitlength}{0.75mm} \thicklines
\begin{picture}(80,15)

\put(0,0){\circle*{2.5}} \put(20,0){\circle*{2.5}}
\put(40,0){\circle*{2.5}} \put(60,0){\circle*{2.5}}
\put(80,0){\circle*{2.5}}

\put(0,0){\line(1,0){80}}

\put(0,5){\makebox(0,0){$2$}} \put(20,5){\makebox(0,0){$1$}}
\put(40,5){\makebox(0,0){$3$}} \put(60,5){\makebox(0,0){$4$}}
\put(80,5){\makebox(0,0){$5$}}

\end{picture}
\end{center}


\noindent but

\begin{center}
\setlength{\unitlength}{0.75mm} \thicklines
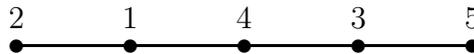
\begin{picture}(80,15)

\put(0,0){\circle*{2.5}} \put(20,0){\circle*{2.5}}
\put(40,0){\circle*{2.5}} \put(60,0){\circle*{2.5}}
\put(80,0){\circle*{2.5}}

\put(0,0){\line(1,0){80}}

\put(0,5){\makebox(0,0){$2$}} \put(20,5){\makebox(0,0){$1$}}
\put(40,5){\makebox(0,0){$4$}} \put(60,5){\makebox(0,0){$3$}}
\put(80,5){\makebox(0,0){$5$}}

\end{picture}
\end{center}

\hspace{1em}

\noindent is not.

A Mathematica routine that computes all the possible $\mu$-labelings
for an order $n$ path is given in Listing \ref{list2}. For example,
the distinct $\mu$-labelings for a path with for 4 vertices are
$$
\begin{array}{cccl}
 \{1,2,3,4\}  & \{1,2,4,3\}  & \{1,4,2,3\}  & \{1,4,3,2\}  \\
 \{2,1,3,4\}  & \{2,3,1,4\}  & \{2,1,4,3\} & \{3,2,1,4\}\,.
\end{array}
$$

\begin{lstlisting}[float=htb,label=list2,caption=Mathematica routine to test and compute $\mu$-labelings]

qPermutations[n_] := Flatten[Table[Flatten[{i, #, j}] & /@
  Permutations[Complement[Range[n], {i, j}]],
    {i, n}, {j, n, i + 1, -1}], 2]

testPermutation[perm_] :=
  Module[{pair1, pair2, i, j, k, l},
    Catch[
      For[p1 = 1, p1 <= Length@perm - 3, ++p1,
        pair1 = Sort@perm[[{p1, p1 + 1}]];
        For[p2 = p1 + 2, p2 < Length@perm, ++p2,
          pair2 = Sort@perm[[{p2, p2 + 1}]];
          If[First@pair1 < First@pair2,
            {i, j} = pair1; {k, l} = pair2,
            {i, j} = pair2; {k, l} = pair1];
          If[ i > k || (j > k && l > j), Throw[False]];
        ]
      ];
      Throw[True]
    ]
  ]

labelings[n_] := Select[qPermutations[n], testPermutation]

\end{lstlisting}

Table \ref{table1} presents the number of distinct $\mu$-labelings
for paths of order up to $11$, which were computed with the same
routines.

\begin{table}[htb]
\caption{\label{table1}Total number of $\mu$-labelings}
\begin{tabular}{c|c}
order & \#labelings \\ \hline
2 & 1\\
3 & 3\\
4 & 8\\
5 & 20\\
6 & 48\\
\end{tabular} \hspace{1cm}
\begin{tabular}{c|c}
order & \#labelings \\ \hline
7 & 112\\
8 & 256\\
9 & 576\\
10 & 1280\\
11 & 2816\\
\end{tabular}
\end{table}

Remarkably, this exhaustive enumeration leads us exactly to the
integer sequence A001792 of the The On-Line Encyclopedia of Integer
Sequences \cite{A001792}. This sequence has many different
interpretations. Originally, we will find it in \cite{JMCP1946} (cf.
also \cite[Table 22.3]{AS1964}) in the absolute value of the
coefficients of $x^n$ for the Chebyshev polynomials of the first
kind $T_{n+2}$. The most simple formula is perhaps $(n+2)2^{n-1}$,
for each positive integer $n$. This sequence emerges also from the
Bernoulli's triangle rows sums \cite{Bernoulli,A008949}.
Nonetheless, we can recent find a vast number of interesting
interpretations for this sequence. Namely, it is the determinant of
the square matrix with 3's on the diagonal and 1's elsewhere, or the
absolute value of the determinant of the Toeplitz matrix with first
row containing the first $n$ integers \cite{A001792}.

Regarding other Mathematica routines, one can find several collected
in \cite{A001792} such as:

\

\texttt{\small matrix[n\_Integer /; n >= 1] := Table[Abs[p - q] + 1, {q,
n}, {p, n}];}

\texttt{\small a[n\_Integer /; n >= 1] := Abs[Det[matrix[n]]]}

\hspace{0.5em}

or

\hspace{0.5em}

\texttt{\small g[n\_, m\_, r\_] := Binomial[n - 1, r - 1] Binomial[m + 1,
r] r;}

\hspace{0.5em}

or

\hspace{0.5em}

\texttt{\small Table[1 + Sum[g[n, k - n, r], {r, 1, k}, {n, 1, k - 1}], {k,
1, 29}]}

\hspace{0.5em}

or

\hspace{0.5em}

\texttt{\small LinearRecurrence[{4, -4}, {1, 3}, 40]}

\hspace{0.5em}

or even

\hspace{0.5em}

\texttt{\small CoefficientList[Series[(1 - x) / (1 - 2 x)\^{}2, {x, 0,
40}], x]}

\hspace{0.5em}

One interesting open question is formally prove that this new
labeling leads to the integer sequence A001792.

\end{document}